\documentstyle[psfig,11pt]{amsart}

\title{On the $K$-theory of elliptic curves}
\author{Kevin P. Knudson}\thanks{Supported by an NSF Postdoctoral
Fellowship, grant no. DMS--9627503} 
\address{Department of Mathematics, Northwestern University, Evanston,
IL  60208}
\email{knudson@@math.nwu.edu}
\date{July 13, 1998, revised September 22, 1998}

\newtheorem{theorem}{Theorem}[section]
\newtheorem{prop}[theorem]{Proposition}

\newtheorem{cor}[theorem]{Corollary}

\newcommand{\X}{{\mathcal X}}
\newcommand{\D}{{\mathcal D}}
\newcommand{\cL}{{\mathcal L}}
\newcommand{\ebar}{\overline{E}}
\newcommand{\knrq}{K_n(R)_{{\mathbb Q}}}
\newcommand{\knaq}{K_n(A)_{{\mathbb Q}}}
\newcommand{\kaq}{K_2(A)_{{\mathbb Q}}}
\newcommand{\knkq}{K_n(k)_{{\mathbb Q}}}
\newcommand{\kkq}{K_2(k)_{{\mathbb Q}}}
\newcommand{\knfq}{K_n(F)_{{\mathbb Q}}}
\newcommand{\kfq}{K_2(F)_{{\mathbb Q}}}
\newcommand{\kneq}{K_n(\overline{E})_{{\mathbb Q}}}
\newcommand{\keq}{K_2(\overline{E})_{{\mathbb Q}}}
\newcommand{\zz}{{\mathbb Z}}
\newcommand{\zq}{{\mathbb Q}}
\newcommand{\uk}{k^\times}
\newcommand{\ukw}{k(\omega)^\times}
\newcommand{\pgla}{PGL_2(A)}
\newcommand{\gla}{GL_2(A)}
\newcommand{\pglk}{PGL_2(k)}
\newcommand{\glk}{GL_2(k)}
\newcommand{\pglf}{PGL_2(F)}
\newcommand{\glf}{GL_2(F)}
\newcommand{\lra}{\longrightarrow}
\newcommand{\ra}{\rightarrow}
\newcommand{\bop}{\bigoplus}

\begin{document}

\begin{abstract} Let $A$ be the coordinate ring of an affine elliptic
curve (over an infinite field $k$) of the form $X-\{p\}$, where 
$X$ is projective and $p$ is a closed point on $X$.  Denote by
$F$ the function field of $X$.  We show that the image of
$H_\bullet(\gla,\zz)$ in $H_\bullet(\glf,\zz)$ coincides with
the image of $H_\bullet(\glk,\zz)$.  As a consequence, we
obtain numerous results about the $K$-theory of $A$ and $X$.
For example, if $k$ is a number field, we show that $r_2(K_2(A)
\otimes \zq) = 0$, where $r_m$ denotes the $m$th level of the
rank filtration.
\end{abstract}
\maketitle

\section{Introduction}\label{intro}

Computing the $K$-theory of a scheme $X$ is a very difficult task.  Even
the simplest case $X=\textrm{Spec}\, k$, where $k$ is a field, is not
completely solved, although a great deal is known.  The next case to 
consider is when $X$ is a curve over $k$, and it is here that the
complexity grows rapidly.  Some curves of genus zero present no real
difficulty thanks to the fundamental theorem:  $K_i(R[t,t^{-1}])=
K_i(R) \oplus K_{i-1}(R)$ for $R$ regular.
The $K$-theory of elliptic curves, on the other hand, has proved to be much 
more elusive.

A great deal of recent work has focused on the construction of specific
elements in the $K$-theory of elliptic curves, particularly in the second
group $K_2$.  This program goes back to the work of S.~Bloch \cite{bloch},
who constructed a regulator map on $K_2$ and used it to find nontrivial
elements.  A.~Beilinson \cite{beil} generalized this construction and made
a number of conjectures relating the dimension of $K_2 \otimes \zq$ with
the values of $L$-functions on the curve.  More recently, Goncharov--Levin 
\cite{gonch}, Rolshausen--Schappacher \cite{rs}, and Wildeshaus \cite{wild}
have made further progress in this area.

In this paper we consider the following situation.  Let $E$ be an
affine elliptic curve defined by the Weierstrass equation $F(x,y)=0$, where
$$F(x,y) = y^2 + a_1xy + a_3y - x^3 - a_2x^2 - a_4x - a_6.$$  Here, the
$a_i$ lie in an infinite field $k$.  Denote by $\ebar$ the projective
curve $E \cup \{\infty\}$ and by $F$ the function field of $\ebar$. 
Denote by $A$ the affine coordinate ring of $E$; it is a Dedekind domain
with field of fractions $F$.  We have $A^\times = k^\times$.

Consider the obvious embedding $i:\gla \lra \glf$.  The main result of
this paper is the following.

\begin{theorem}\label{mainthm}  The image of the map
$$i_*:H_\bullet(\gla,\zz)\lra H_\bullet(\glf,\zz)$$
coincides with the image of
$$(i|_{\glk})_*:H_\bullet(\glk,\zz)\lra H_\bullet(\glf,\zz).$$
\end{theorem}

This is a consequence of an explicit computation of the homology of
$\pgla$ due to the author \cite{knudson} (recalled in Section 
\ref{pglhomology} below).  The proof of Theorem \ref{mainthm} is given
in Section \ref{proof}.

\medskip

\noindent {\bf Remark.}  Theorem \ref{mainthm} and its corollaries
in Sections \ref{filt} and \ref{numfld} are valid also for singular
cubic curves $F(x,y)=0$.  We shall point out the necessary modifications
needed to prove this below.

\medskip

From this result we deduce a number of facts about the $K$-theory of
$E$ and $\ebar$.  Recall the {\em rank filtration} of the rational
$K$-theory $K_\bullet(R)_{\zq}:= K_\bullet(R)\otimes_{\zz} \zq$ of
a ring $R$:
$$r_m\knrq = \textrm{im}\{H_n(GL_m(R),\zq) \lra H_n(GL(R),\zq)\}
\cap \knrq.$$

\begin{cor} The image of the map $r_2\knaq \lra r_2\knfq$ coincides with
the image of $r_2\knkq \lra r_2\knfq$.
\end{cor}

In particular, when $n=2$ we see that the image of $r_2\kaq \lra r_2\kfq$
coincides with the image of $\kkq$.

\medskip

\noindent {\bf Remark.}  This corollary is valid for {\em any} field
$k$.  Indeed, if $k$ is finite, then the rational homology
$H_\bullet(\gla,\zq)$ vanishes in positive degrees (as does
$H_\bullet(\glk,\zq)$) from which it follows that $r_2\knaq = 0$.

\medskip

Define a filtration $r_\bullet K_\bullet(\ebar)_{\zq}$ by pulling back
the rank filtration of $K_\bullet(A)_{\zq}$:
$$r_m\kneq := (f^*)^{-1}(r_m\knaq),$$
where $f:E\lra \ebar$ is the inclusion and $f^*:K_\bullet(\ebar) \lra
K_\bullet(E) = K_\bullet(A)$ is the induced map in $K$-theory.  Then
we obviously have the following result.

\begin{cor}  The image of $r_2\kneq \lra r_2\knfq$ coincides with the
image of $r_2\knkq \lra r_2\knfq$.
\end{cor}

We study the filtration $r_\bullet$ in greater detail in Section \ref{filt}.
In Section \ref{numfld} we specialize to the case where $k$ is a number
field.  In this case, we show that $r_2\kaq = 0$.

In the case $n=2$, results of Nesterenko--Suslin \cite{nessus} imply that
$r_3\kaq = \kaq$.  A description of the homology of $PGL_3(A)$ (or $GL_3(A)$)
would provide a great deal of insight into the structure of $\keq$,
especially over a number field.  Such a computation remains elusive,
however.

\medskip

{\small
\noindent {\em Acknowledgements.}
I thank Dick Hain for practically insisting that I study this question.
Andrei Suslin provided valuable insight regarding the proof of Theorem
\ref{mainthm}.  Thanks also go to Mark Walker for answering several
(perhaps silly) questions.}

\section{The Rank Filtration}\label{filt}

The rational $K$-groups of affine schemes admit the rank filtration
mentioned in the introduction.  Since $BGL(R)^+$ is an $H$-space, the
Milnor--Moore Theorem \cite{milmor} implies that the Hurewicz map
$$h:\pi_\bullet(BGL(R)^+)\otimes\zq \lra H_\bullet(GL(R),\zq)$$
is injective with image equal to the primitive elements of the homology.
The rank filtration is the increasing filtration defined by
$$r_m\knrq = \textrm{im}\{H_n(GL_m(R),\zq)\lra H_n(GL(R),\zq)\}
\cap \knrq.$$  
By Theorem 2.7 of \cite{nessus}, the map $H_2(GL_3(A),\zz) \ra H_2(GL(A),\zz)$
is surjective so that $r_3\kaq = \kaq$.  The rank filtration of $\kaq$
then has the form
$$0=r_1\kaq \subseteq r_2\kaq \subseteq r_3\kaq = \kaq$$
(the vanishing of $r_1$ is a consequence of the vanishing of $r_1\kkq$
for infinite fields \cite{nessus}, and the fact that $A^\times = k^\times$).

Define an increasing filtration $r_\bullet$ of $\kneq$ as follows.  Let
$f:E\lra \ebar$ be the canonical inclusion and denote by $f^*$ the
induced map on $K$-theory.  We define $r_m\kneq$ by
$$r_m\kneq = (f^*)^{-1}(r_m\knaq).$$  There is a commutative diagram
\begin{equation}\label{rankdiagram}
\begin{array}{ccc}
r_m\kneq & \lra & r_m\knaq  \\
         &\searrow & \downarrow  \\
         &    &  r_m\knfq.     
\end{array}
\end{equation}

\begin{prop}\label{rkimage} 
The image of $r_2\knaq \lra r_2\knfq$ coincides with
the image of $r_2\knkq \lra r_2\knfq$.
\end{prop}

\begin{pf}  By Theorem \ref{mainthm}, the image of 
$$i_*:H_\bullet(\gla,\zz)\lra H_\bullet(\glf,\zz)$$ coincides with
the image of $(i|_{\glk})_*$.  Consider the commutative diagram
$$\begin{array}{ccc}
H_n(\gla,\zq) & \lra & H_n(GL(A),\zq)   \\
\downarrow    &      & \downarrow       \\
H_n(\glf,\zq) & \lra & H_n(GL(F),\zq).
\end{array}$$
It follows that the image of $H_n(\gla,\zq)$ in $H_n(GL(F),\zq)$
coincides with the image of $H_n(\glk,\zq)$; {\em i.e.,} the image
of $r_2\knaq \ra r_2\knfq$ coincides with the image of $r_2\knkq$.
\end{pf}

\begin{cor}  The image of $r_2\kneq \lra r_2\knfq$ coincides with
the image of $r_2\knkq$.
\end{cor}

\begin{pf}  This follows by considering the diagram (\ref{rankdiagram}).
\end{pf}

\section{The Number Field Case}\label{numfld}

Suppose that the ground field $k$ is a number field.  By localizing
the projective curve at its generic point we obtain the following
exact sequence for $K_2$
$$0 \lra \keq \lra \kfq \stackrel{{\mathcal T}}{\lra}
\bop_P K_1(k(P))_{\zq}$$ where $P$ varies over the closed points of
$\ebar$ and $k(P)$ is the residue field at $P$.  The map ${\mathcal T}$
is the {\em tame symbol} (see, {\em e.g.,} \cite{rs}).

\medskip

\noindent {\bf Remark.}  It is not known for a single curve if
$\keq$ is finite dimensional.  Beilinson has conjectured that the dimension
of this space is related to special values of $L$-functions on 
$\overline{E}$.  This conjecture was modified by Bloch and Grayson 
\cite{blochgrayson} to predict that the dimension
is the number of infinite places of $k$ plus the number of
primes ${\mathfrak p} \subset {\mathcal O}_k$ where $\ebar$ has 
split multiplicative reduction modulo ${\mathfrak p}$.  For a discussion
of this see, for example, \cite{rs}.

\medskip

We also have the localization sequence for $A$:
$$\cdots \ra K_{i+1}(F) \ra \bop_{{\mathfrak p}\; \textrm{maximal}} 
K_i(A/{\mathfrak p}) \ra K_i(A) \ra K_i(F) \ra \cdots .$$  Since
$A/{\mathfrak p}$ is a finite extension of $k$ for all ${\mathfrak p}$,
the groups $K_{2m}(A/{\mathfrak p})$ are torsion.  It follows that we have
an exact sequence
$$0 \lra K_{2m}(A)_{\zq} \lra K_{2m}(F)_{\zq}
 \lra \bop_{\mathfrak p} 
K_{2m-1}(A/{\mathfrak p})_{\zq}.$$

\begin{prop}  If the ground field $k$ is a number field, then
the map $\keq \ra \kaq$ is injective.
\end{prop}

\begin{pf} This follows by considering the commutative diagram
$$\begin{array}{ccccc}
0 & \lra & \kaq & \lra & \kfq \\
  &      &\uparrow &   &  ||  \\
0 & \lra & \keq & \lra & \kfq.
\end{array}$$
\end{pf}

\begin{prop} If $k$ is a number field, then $r_2\kaq = 0 = r_2\keq$.
\end{prop}

\begin{pf}  The map $\kaq \ra \kfq$ is injective.  But by Proposition
\ref{rkimage}, the image of $r_2\kaq$ coincides with the image of
$r_2\kkq = \kkq = 0$.
\end{pf}

As a consequence we see that any nontrivial elements of $\kaq$ (and hence of
$\keq$) must come from $H_2(GL_3(A),\zq)$.  Thus, to prove that $\keq$ is
a finite dimensional vector space, it suffices to show that the 
image of $H_2(GL_3(A),\zq)$
in $H_2(GL_3(F),\zq)=H_2(GL(F),\zq)$ is finite dimensional.

\section{The Homology of $\pgla$}\label{pglhomology}

The remainder of the paper is devoted to the proof of Theorem
\ref{mainthm}.  We begin by recalling the calculation of
$H_\bullet(\pgla,\zz)$ given in \cite{knudson}.  The proof
uses the action of $\pgla$ on a certain Bruhat--Tits tree
$\X$.  

We use the description of $\X$ given by Takahashi \cite{takahashi}.
Recall that $A$ is the coordinate ring of the affine curve $E$ with
function field $F$.  The field $F$ has transcendence degree 1 over
$k$ and is equipped with the discrete valuation at $\infty$, 
$v_{\infty}$.  Denote by ${\mathcal O}_\infty$ the valuation ring
and by $t=x/y$ the uniformizer at $\infty$.  Denote by $\cL$ the
field of Laurent series in $t$ and let $v$ be the valuation on 
$\cL$ defined by $v(\sum\limits_{n\ge n_0} a_nt^n) = n_0$.  The
ring $A$ can be embedded in $\cL$ in such a way that $v(x)=-2$ and
$v(y)=-3$; we identify $A$ with its image in $\cL$.  Note that this
embedding induces an embedding $F\ra \cL$ and that the completion of
$F$ with respect to $v_\infty$ is $\cL$. We therefore have a 
commutative diagram
$$\begin{array}{ccc}
\gla & \lra & \glf  \\
     &\searrow & \downarrow \\
     &         & GL_2(\cL).
\end{array}$$

Let $G=GL_2(\cL)$ and $K=GL_2(k[[t]])$.  Denote by $Z$ the center of
$G$.  The Bruhat--Tits tree $\X$ is defined as follows.  The vertex set
of $\X$ is the set of cosets $G/KZ$.  Two cosets $g_1KZ$ and $g_2KZ$
are adjacent if 
$$g_1^{-1}g_2 = \left(\begin{array}{cc}
                        t  &  b  \\ 
                        0  &  1
                      \end{array}\right) \quad \textrm{or} \quad
                 \left(\begin{array}{cc}
                       t^{-1} & 0 \\
                         0    & 1
                      \end{array}\right)  \quad \textrm{modulo}\, KZ$$
for some $b\in k$.  The graph $\X$ is a tree \cite{serre}, p.~70.
Note that $\gla$ acts on $\X$ without inversion and that the center of
$\gla$ (which is equal to $\uk$) acts trivially on $\X$.  It follows
that the quotient $\pgla\backslash\X$ is defined.  We describe a
fundamental domain $\D \subset \X$ for the action ({\em i.e.,}
$\D \cong \pgla\backslash\X$).

If $f_1,f_2\in \cL$, denote by $\phi(f_1,f_2)$ the vertex
$\left(\begin{array}{cc}
    f_1 & f_2 \\
     0  &  1
       \end{array}\right) KZ$.  Denote by $F_x(l,m)$ and $F_y(l,m)$
the partial derivatives at $(l,m)$ of the Weierstrass equation
$F(x,y)$.  Define two sets $E_1$ and $E_2$ as follows:
$$E_1 = \{(l,m): F(l,m)=0 \;\textrm{and}\; F_y(l,m)=0\} \cup \{\infty\}$$
and
$$E_2 = \{(l,m): F(l,m)=0 \;\textrm{and}\; F_y(l,m)\ne 0\}.$$
Observe that $\ebar$ = $E_1 \cup E_2$.  Define vertices of $\X$ by
\begin{eqnarray*}
o   &  =   &   \phi(t,t^{-1});    \\
v(l)&  =   & {\begin{cases}
                \phi(t^2,t^{-1}+lt)  &  \textrm{if}\; l\in k \\
                \phi(1,t^{-1})       &  \textrm{if}\; l=\infty;
              \end{cases}}    \\
c(p,n)& =  & {\begin{cases}
                \phi(t^{n+2},\frac{y-m}{x-l})  &  \textrm{if}\; p=(l,m)\in E \\
                \phi(t^{-n},0)          & \textrm{if}\; p=\infty;
              \end{cases}}    \\
e(p)  & =  & {\begin{cases}
               \phi(t^4,\frac{y-m}{x-l} + \frac{F_x(l,m)}{y-m})  &
                                           \textrm{if}\; p=(l,m)\in E_1 \\
               \phi(1,0)                   &  \textrm{if}\; p=\infty.
               \end{cases}}
\end{eqnarray*}

We are now ready to describe the subgraph $\D$.  For each $l \in k \cup
\{\infty\}$, the vertex $v(l)$ is adjacent to $o$.  Denote by
$\D(l)$ the connected component of $\D-\{o\}$ which contains $v(l)$.
The $\D(l)$ fall into three types.

(1) Suppose $F(x,y)=0$ has no rational solution with $x=l$.  Then 
$\D(l)$ consists only of $v(l)$ (see Figure 1).

\begin{figure}
\centerline{\psfig{figure=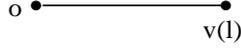,angle=270,width=1.25in}}
\caption{$F(l,y)=0$ has no rational solutions}
\end{figure}

(2) Suppose $l=\infty$ or $F(x,y)=0$ has a unique rational solution
with $x=l$.  Let $p$ be the point at infinity of $E$ or the rational
point corresponding to the solution.  Note that $p$ is a point of
order $2$.  Then $\D(l)$ consists of an infinite path $c(p,1),c(p,2),
\dots$ and an extra vertex $e(p)$ (see Figure 2).

\begin{figure}
\centerline{\psfig{figure=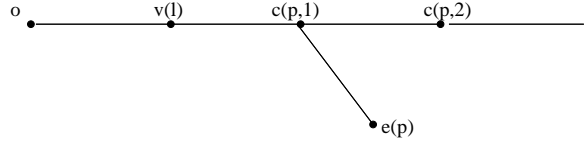,angle=270,width=3in}}
\caption{$F(l,y)=0$ has a unique rational solution}
\end{figure}

(3)  Suppose $F(x,y)=0$ has two different solutions such that $x=l$.
Let $p,q$ be the corresponding points on $E$.  Then $\D(l)$ consists 
of two infinite paths $c(p,1),c(p,2),\dots$ and
$c(q,1),c(q,2),\dots$ (see Figure 3).

\begin{figure}
\centerline{\psfig{figure=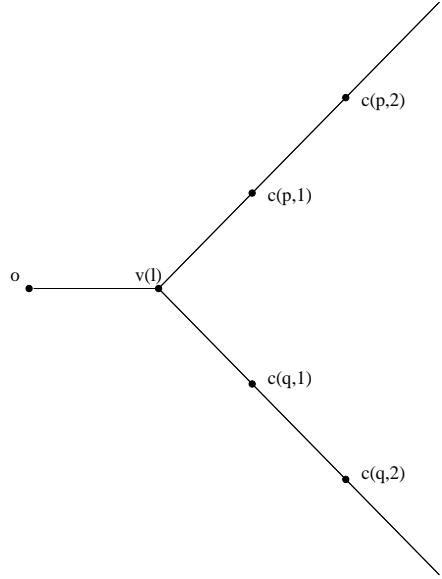,angle=270,height=3in,width=2.25in}}
\caption{$F(l,y)=0$ has two distinct solutions}
\end{figure}

The infinite path $c(p,1),c(p,2),\dots$ is called a {\em cusp}.  Note
that there is a one-to-one correspondence between cusps and the
rational points of $\overline{E}$.

\begin{theorem}[Takahashi]  The graph $\D$ is a fundamental domain
for the action of $\gla$ on $\X$ (and hence is also a fundamental
domain for the action of $\pgla$). \hfill $\qed$
\end{theorem}

\medskip

\noindent {\bf Remark.}  The theorem is true also for singular curves
$C$ given by
$F(x,y)=0$ with the following modification.  If the curve is singular
at $p=(l,m)$, then the vertex $e(p)$ is the same as $c(p,2)$.  In
this case, then, the tree $\D(l)$ consists only of the cusp
$c(p,1),c(p,2),\dots$.  The proofs of the following results for $C$
then go through unchanged except that the summands in the homology
decomposition of $H_\bullet(PGL_2(k[C]),\zz)$ corresponding to
singular points are $H_\bullet(\uk,\zz)$ instead of $H_\bullet(\pglk,\zz)$.

\medskip

Since $\X$ is contractible, we have a spectral sequence with $E^1$--term
$$E_{p,q}^1 = \bigoplus_{\sigma^{(p)} \subset \D}
H_q(\Gamma_{\sigma},\zz) \Longrightarrow H_{p+q}(\pgla,\zz)$$
where $\Gamma_{\sigma}$ is the stabilizer of the $p$--simplex $\sigma$
in $\pgla$.  We shall discuss the stabilizers in detail in the next
section.  For the purpose of computing homology, the next result is 
sufficient (see \cite{takahashi}, Theorem 5).  If $F(l,y)=0$ has
no rational solution, denote by $k(\omega)$ the quadratic extension
of $k$ in which $F(l,\omega)=0$.

\begin{prop}\label{stabhom} Up to isomorphism, the stabilizers $\Gamma_\sigma$
are as follows:
\begin{eqnarray*}
\Gamma_o   &   =   &   \{1\}     \\
\Gamma_{v(l)} & \cong & {\begin{cases}
                       \ukw/\uk    &   \text{\em in case (1)} \\
                         k         &   \text{\em in case (2)} \\
                         \uk       &   \text{\em in case (3)}
                     \end{cases}}  \\
\Gamma_{c(p,n)}& \cong & {\biggl\{ \left(\begin{array}{cc}
                                            p  &  v \\
                                            0  &  q
                                         \end{array}\right): p,q \in \uk,
                                             v\in k^n\biggr\}/\uk}   \\
\Gamma_{e(p)} & \cong & \pglk.
\end{eqnarray*}
The stabilizer of an edge is the intersection of its vertex stabilizers
(one of which is contained in the other). \hfill $\qed$
\end{prop}

By Theorem 1.11 of \cite{nessus}, the inclusion of the diagonal subgroup
into $\Gamma_{c(p,n)}$ induces an isomorphism in homology.  This
leads to the proof of the following, which is the main result of
\cite{knudson}.

\begin{theorem}\label{pglhom}
For all $i \ge 1$,
\begin{eqnarray*}
H_i(\pgla,\zz) & \cong & {\bop\begin{Sb}
                           l\in k\cup\{\infty\}\\
                         F(l,y)=0\,\text{\em has unique sol.}
                         \end{Sb}  
                      H_i(\pglk,\zz)}   \\
               &     &  {\oplus
                       \bop\begin{Sb}
                            l\in k \\
                          F(l,y)=0\,\text{\em has two sol.}
                           \end{Sb}
                      H_i(\uk,\zz)}  \\
                &   &  {\oplus \bop\begin{Sb}
                            l\in k \\
                          F(l,y)=0\,\text{\em has no sol.}
                              \end{Sb}
                          H_i(\ukw/\uk,\zz).}   \hfill \qed
\end{eqnarray*}
\end{theorem}

\medskip

\noindent {\bf Remark.}  This theorem holds also in degrees $\le 2$ if
$k$ is a finite field with at least $4$ elements.  For in this case, the
inclusion of the diagonal subgroup into $\Gamma_{c(p,n)}$ induces a
homology isomorphism in degrees $\le 2$; see \cite{rosenberg}, p.~204.

\medskip

The isomorphism is induced by the inclusion of the various 
$\Gamma_{v(l)}$ and $\Gamma_{e(p)}$.  In the next section, we shall
compute the image of the map
$$H_\bullet(\pgla,\zz)\lra H_\bullet(\pglf,\zz).$$

\section{The Map $H_\bullet(\pgla,\zz)\ra H_\bullet(\pglf,\zz)$}
\label{map}

To compute the image of $H_\bullet(\pgla,\zz)$ in $H_\bullet(\pglf,\zz)$,
we must examine the various $\Gamma_v$ in greater detail.  If
$p=\infty$, then the stabilizer $\Gamma_{e(\infty)}$ is the subgroup
$\pglk$ of $\pgla$.  Hence, under the map $j:\pgla\ra\pglf$, 
$\Gamma_{e(\infty)}$ maps to $\pglk \subset \pglf$.

The other stabilizers for $l\ne \infty$ are {\em not} subgroups of
$\pglk$, although they are isomorphic to such.  We have the following
result.

\begin{theorem}\label{conj}  For each $l\in k$, the stabilizers
$\Gamma_{v(l)}$ and $\Gamma_{e(p)}$ ($p=(l,m)$) are conjugate in
$\pglf$ to subgroups of $\pglk$.
\end{theorem}

\begin{cor}\label{image} The image of $j_*:H_\bullet(\pgla,\zz)
\ra H_\bullet(\pglf,\zz)$ coincides with the image of
$H_\bullet(\pglk,\zz)$.
\end{cor}

\begin{pf}  It is well-known (see \cite{brown}, p.~48) that conjugation
induces the identity on homology.  It follows that if $H_1,H_2$ are
conjugate subgroups of a group $G$, then the images of $H_\bullet(H_i,\zz)
\ra H_\bullet(G,\zz)$ coincide.  Since each stabilizer which appears
in the homology decomposition of $\pgla$ is conjugate in $\pglf$ to
a subgroup of $\pglk$, the result follows.
\end{pf}

\noindent {\em Proof of Theorem} \ref{conj}.  To keep the notation as
simple as possible, we only prove the case $\Gamma_{v(0)}$ and in
the case $F(0,0)=0=F_y(0,0)$, $\Gamma_{e(0,0)}$.  All other cases
are similar (but notationally more complex).  For $r_1,\dots ,r_4 \in k$
define
$$M_2(r_1,r_2) = \left(\begin{array}{cc}
                 r_2y + r_1   &   -r_2\biggl(\frac{y^2+a_3y-a_6}{x}\biggr) \\
                 r_2x         & -r_2y-a_3r_2+r_1
                       \end{array}\right)$$
and
\begin{eqnarray*}
\lefteqn{M_4(r_1,r_2,r_3,r_4) =} \\
 & & { \left(\begin{array}{c|c}
            r_4xy+r_3(x^2+a_2x+a_4)  & {-r_4y^2 -r_3y(x+a_2)+a_4r_4(x+a_2)} \\
                    + r_2y + r_1     &   {-r_2(x^2+a_2x+a_4-a_1y)} \\ \hline
           r_4x^2 + r_3(y+a_1x) & {-r_4xy-r_3(x^2+a_2x+a_4)}\\
               + r_2x + a_4r_4   &   {-r_2y +a_1a_4r_4+a_4r_3+r_1}
                             \end{array}\right).}
\end{eqnarray*}
According to Proposition 9 of \cite{takahashi}, the stabilizer of
$v(0)$ in $\gla$ is
$$\tilde{\Gamma}_{v(0)} = \{M_2(r_1,r_2): r_1(-a_3r_2+r_1)-a_6r_2^2 \ne 0\},$$
and of $e(0,0)$ is
\begin{eqnarray*}
\lefteqn{\tilde{\Gamma}_{e(0,0)} =} \\
 &  & { \{M_4(r_i): r_1(a_4r_3 + r_1)
                             + (-a_2a_4r_4+a_1a_4r_3+a_4r_2a_1r_1)a_4r_4 
                                 \ne 0 \}.}
\end{eqnarray*}
Consider the following identity:
$$\left(\begin{array}{cr}
               x  & -y  \\
               0  &  1
        \end{array}\right) M_2(r_1,r_2) \left(\begin{array}{cc}
                                             1/x   &  y/x  \\
                                              0    &   1
                                            \end{array}\right)
= \left(\begin{array}{cc}
         r_1  &  a_6r_2  \\
          r_2  & r_1-a_3r_2
       \end{array}\right) = N_2(r_1,r_2).$$
It follows that
$$\left(\begin{array}{cr}
               x  & -y  \\
               0  &  1
        \end{array}\right)  \tilde{\Gamma}_{v(0)}\left(\begin{array}{cc}
                                             1/x   &  y/x  \\
                                              0    &   1
                                            \end{array}\right)
= \{N_2(r_1,r_2):\det N_2(r_1,r_2)\ne 0\} = \tilde{\Gamma}.$$
Note that the subgroup $\tilde{\Gamma}$ lies in $\glk$ and that 
$g=\left(\begin{array}{cr}
          x  &  -y  \\ 0  &  1
         \end{array}\right)$ is an element of $\glf$.  It follows that
$g\Gamma_{v(0)}g^{-1} = \tilde{\Gamma}/\uk \subset \pglk$ 
inside $\pglf$.  Moreover,
we can demonstrate the isomorphism of Proposition \ref{stabhom} as follows.
If $F(0,y)=0$ has no rational solutions, then define a map $\tilde{\Gamma}
\ra \ukw$ by $N_2(r_1,r_2) \mapsto r_1 + r_2\omega$.  One checks easily
that this is an isomorphism.  If $F(0,y)=0$ has two solutions, say
$u,v\in k$, then it is easy to see that
$$\left(\begin{array}{cc}
       u  &   uv \\  \frac{-1}{u(u-v)} & \frac{-1}{u-v}
        \end{array}\right)\tilde{\Gamma} \left(\begin{array}{cc}
       u  &   uv \\  \frac{-1}{u(u-v)} & \frac{-1}{u-v}
        \end{array}\right)^{-1} = D(k)$$
where $D(k) \subset \glk$ is the subgroup of diagonal matrices.  Finally,
if $F(0,y)=0$ has a unique solution, say $u\in \uk$, then
$$\left(\begin{array}{cr}
          0  &  -1  \\  1  &  u  
        \end{array}\right) \tilde{\Gamma} \left(\begin{array}{cr}
          0  &  -1  \\  1  &  u  
        \end{array}\right)^{-1} = B(k)$$
where $B(k)$ is the upper triangular subgroup of $\glk$.

For the group $\Gamma_{e(0,0)}$ we have
$$\left(\begin{array}{cc}
       x  &  -y  \\  \frac{y}{a_4} & 1 - \frac{y^2}{a_4x}
        \end{array}\right) \tilde{\Gamma}_{e(0,0)}                                            \left(\begin{array}{cc}
       x  &  -y  \\  \frac{y}{a_4} & 1 - \frac{y^2}{a_4x}
        \end{array}\right)^{-1} = \glk$$
from which it follows that $\Gamma_{e(0,0)}$ is conjugate to $\pglk$
inside $\pglf$.  (Note that since $\ebar$ is smooth, $a_4\ne 0$.)
\hfill $\qed$

\section{Proof of Theorem \ref{mainthm}}\label{proof}

We now prove Theorem \ref{mainthm}. Corollary \ref{image}
shows that $H_\bullet(\pgla,\zz)$ has image equal to the image
of $H_\bullet(\pglk,\zz)$ in $H_\bullet(\pglf,\zz)$.  Consider the
following commutative diagram
$$\begin{array}{ccccccccc}
1 & \ra & \uk & \ra & \gla & \ra & \pgla & \ra & 1 \\
  & & \downarrow &  & \downarrow & & \downarrow & & \\
1 & \ra & F^\times & \ra & \glf & \ra & \pglf & \ra & 1
\end{array}$$
and the induced map of Hochschild--Serre spectral sequences
$$\begin{array}{ccccc}
E_{p,q}^2(A) & = & H_p(\pgla,H_q(\uk)) & \Rightarrow & H_{p+q}(\gla,\zz) \\
\downarrow   &   &  \downarrow         &             &  \downarrow  \\
E_{p,q}^2(F) & = & H_p(\pglf,H_q(F^\times)) & \Rightarrow & H_{p+q}(\glf,\zz).
\end{array}$$
Since the extensions are central, the groups $H_q(\uk)$ (resp.
$H_q(F^\times)$) are trivial $\pgla$ (resp. $\pglf$) modules.  Hence
we have the following commutative diagram of universal coefficient
sequences
$$\begin{array}{ccccc}
\scriptstyle{H_p(\pgla)\otimes H_q(\uk)} & \scriptstyle{\ra} &
\scriptstyle{H_p(\pgla,H_q(\uk))} & \scriptstyle{\ra} & 
\scriptstyle{\textrm{Tor}_1^{\zz}(H_{p-1}(\pgla),H_q(\uk))} \\
\scriptstyle{\downarrow} & & \scriptstyle{\downarrow} & & 
\scriptstyle{\downarrow}  \\
\scriptstyle{H_p(\pglf)\otimes H_q(F^\times)} & \scriptstyle{\ra} & 
\scriptstyle{H_p(\pglf,H_q(F^\times))} & \scriptstyle{\ra} & 
\scriptstyle{\textrm{Tor}_1^{\zz}(H_{p-1}(\pglf),H_q(F^\times)).}
\end{array}$$
By Corollary \ref{image}, we see that the image of $E_{p,q}^2(A) \ra
E_{p,q}^2(F)$ coincides with the image of $E_{p,q}^2(k) \ra E_{p,q}^2(F)$.
It follows that the same is true of the $E^{\infty}$ terms:
$$\textrm{im}\{E_{p,q}^{\infty}(A)\ra E_{p,q}^{\infty}(F)\}
=\textrm{im}\{E_{p,q}^{\infty}(k)\ra E_{p,q}^{\infty}(F)\}.$$
Thus, the image of $H_\bullet(\gla,\zz)$ in $H_\bullet(\glf,\zz)$
coincides with the image of $H_\bullet(\glk,\zz)$.

\end{document}